\input amstex
\documentstyle{amsppt}
\magnification 1200
\NoBlackBoxes
\NoRunningHeads
%\operatorname{cotg}
%\operatorname\c\cotan
\topmatter
\title On martingale approximation of adapted processes
\endtitle
\author Herv\'e Queff\'elec and Dalibor Voln\'y \endauthor
\abstract
We show that the existence of a martingale approximation of a stationary process depends on the choice of the filtration. There 
exists a stationary linear process which
has a martingale approximation with respect to the natural filtration, but no approximation with respect to a larger filtration
with respect to wich it is adapted and regular. There exists a stationary process adapted, regular, and having a martingale
approximation with respect to a given filtration but not (regular and having a martingale approximation) with respect
to the natural filtration.
\endabstract
\endtopmatter
\document

\subheading{1. Introduction}
Let $(\Omega,\Cal A,\mu)$ be a probability space with a bimeasurable and measure-preserving transformation $T$.
For simplicity we shall suppose that $\mu$ is ergodic, i.e\. for $A$ measurable, $A = T^{-1}A$ implies that 
$\mu(A)=0$ or $\mu(A)=1$. 
For any measurable function $f$, $(f\circ T^i)$ is a strictly stationary process. Let $(\Cal F_i)$ be an increasing 
filtration of sub-$\sigma$-algebras of $\Cal A$ such that $\Cal F_i = T^{-i}\Cal F_0$. By $\Cal F_{-\infty}$ we denote 
the intersection of all $\Cal F_i$ and by $\Cal F_{\infty}$ the $\sigma$-algebra generated by $\Cal F_i$, $i\in\Bbb Z$.

We say that for the process $(f\circ T^i)$ there is a {\it martingale approximation w.r.t\. $(\Cal F_i)$} if there 
exists $m\in L^2(\Cal F_0) \ominus L^2 (\Cal F_{-1})$ such that 
$\frac{1}{\sqrt n} \|S_n(f-m)\|_2 \to 0$ for $n\to\infty$ where $S_{n}(g)=\sum_{i=0}^{n-1}g\circ T^{i}$.

This condition is equivalent to Gordin's 1969 condition (cf\. \cite{V93}). The martingale approximation, 
since Gordin's 1969 paper \cite{Go}, has been a powerful tool in proving central limit theorems for stationary processes.
Most results from 70's appear in the monography \cite{Ha-He}, for several more recent contributions we can quote
e.g\. \cite{DeM}, \cite{M-Wo} \cite{Wu-Wo} \cite{Pe-U} \cite{Wu07} \cite{Z-Wo} \cite{V93}. 
The method has been particularly suitable for the study of 
processes of $X_k = g(\dots, e_{-1},e_0)$, where 
$(e_i)$ is an iid sequence; this class includes the (functionals of) stationary linear processes and an interesting study
of limit theorems for this class of processes have been done by W. B. Wu in \cite{Wu05}.

Here we shall deal with the problem of the dependence of the approximation on the choice of the filtration.
The problem was posed by Gordin in his 1968 Ph.D thesis but seemingly, except \cite{V09}, no results have been published.
It has only been known as a folklore theorem, that there exist processes which are deterministic, i.e\. measurable w.r.t\.
the natural filtration, but having a martingale approximation w.r.t\. another filtration (cf\. \cite{V09}).
In \cite{V09} a stationary linear process $(X_i)$ with independent innovations $e_i$ has been found such that
$(X_i)$ is a martingale difference sequence (a sequence of iid, in fact) but does not admit any martingale approximation
w.r.t\. the filtration $(\Cal F_i)$ given by the process of innovations. The example, however, is not adapted to the 
filtration $(\Cal F_i)$.
Here we shall treat a more difficult situation, when the process is adapted. 

We shall say that the process $(f\circ T^i)$ is regular w.r.t\. a filtration $(\Cal F_i)$ if $f$ is 
$\Cal F_{\infty}$-measurable and $E(f|\Cal F_{-\infty}) = 0$.
Let us suppose that the process $(f\circ T^i)$ is adapted w.r.t\. $(\Cal F_i)$, i.e\. $f$ is $\Cal F_0$-measurable.
Denote by $\Cal G_j = \sigma\{f\circ T^i| i\leq j\}$ the natural filtration. Recall that by a stationary linear process
$(X_k)$ with innovations $e_i$ we understand the process of $X_k = \sum_{i\in \Bbb Z} a_ie_{k-i}$ where 
$\sum_{i\in \Bbb Z} a_i^2 <\infty$.

\proclaim{Proposition 1} There exists an adapted stationary linear process $(X_k)$ such
that the $X_k$'s are iid and if  $(\Cal F_i)$ is the filtration given by the innovations, there is a $c>0$ such 
that for every $m\in L^2(\Cal F_0) \ominus L^2 (\Cal F_{-1})$
$$
  \sup_{m\in L^2(\Cal F_0) \ominus L^2 (\Cal F_{-1})} \limsup_{n\to\infty} \frac1{\sqrt n} \|\sum_{i=0}^{n-1}(X_i-m)\|_2 
  >c. \tag{1}
$$
\endproclaim

The proposition shows that $(X_k)$ is also a stationary linear process with another process $(g_i=X_i)_i$ of independent
innovations, adapted to the natural filtration $(\Cal G_i)$, and for $(\Cal G_i)$ we can get much better approximation 
properties. We can easily change the process  $(X_n)$ (e.g\. by adding a coboundary) so that it is not an iid sequence 
(nor a sequence of martingale differences), having no martingale approximation w.r.t\. $(\Cal F_i)$ and having a
martingale approximation w.r.t\. the natural filtration.
\medskip

One can guess that the natural filtration will bring the best approximation. This, however, is not the case.

\proclaim{Proposition 2} 
There exists a process $(f\circ T^i)$ and a filtration $(\Cal F_j)$ such that
\roster
\item"(i)" $(f\circ T^i)$ is not regular w.r.t\. the natural filtration $(\Cal G_i)$,
\item"(ii)" there is no martingale approximation w.r.t\. the natural filtration $(\Cal G_i)$,
\item"(iii)" there is a martingale approximation w.r.t\. the filtration $(\Cal F_i)$.
\endroster
\endproclaim	

\underbar{Remark 1.} The example proving Proposition 2 is such that $g = E(f|\Cal G_{-\infty})$ is a function different 
both from $f$ and from $0$, such that $(g\circ T^i)$ is a nontrivial iid sequence.

In general, we can take the natural filtration $(\Cal G_i^{(1)})_i$ of the process $(E(f|\Cal F_{-\infty})\circ T^i)$, then
of the process $(E(f|\Cal G_{-\infty}^{(1)})\circ T^i)$, ... and get an ordered set of $\sigma$-algebras 
$\Cal G_{-\infty}^{(\alpha)}$, $\alpha < \alpha_0$, where $\alpha_0$ is a countable ordinal number, cf\. \cite{V85} 
(in the exemple proving Proposition 2, the process stops already after the second step).
For such ``cascades'' of filtrations we can have martingale approximations and limit theorems similarly as in the case of 
single filtrations (see \cite{V92}).

\underbar{Remark 2.} There exists a ``deterministic process'' $(f\circ T^i)$ such that $f\circ T^i$ are all measurable
w.r.t\. the past $\sigma$-algebra $\sigma\{f\circ T^i | i\leq 0\}$ hence there is no martingale approximation w.r.t\. the
natural filtration but there is a martingale approximation w.r.t\. another filtration (cf\. \cite{V09}). 

It remains an open question whether there exists a process $(f\circ T^i)$ regular w.r.t\. the natural filtration, without 
martingale 
approximation w.r.t\. the natural filtration, but with a martingale approximation w.r.t\. another filtration $(\Cal F_i)$.

What can be shown for the moment is a process having a martingale approximation both w.r.t\. the natural filtration and
w.r.t\. another filtration $(\Cal F_i)$, and the rate of approximation for $(\Cal F_i)$ is much faster.

\proclaim{Proposition 3} 
For any sequences of $a_n>0, b_n>0$, $(a_n)$ non decreasing, $b_n \searrow 0$, there exists a process $(f\circ T^i)$ and a filtration 
$(\Cal F_j)$, $m'\in L^2(\Cal F_0)\ominus L^2(\Cal F_{-1})$, such that $\|S_n(f-m')\|_2 \leq a_n$ for all $n\geq 1$, 
while for the natural filtration $\Cal G_j = \sigma(f\circ T^i : i\leq j)$ and any $m''\in L^2(\Cal G_0)\ominus 
L^2(\Cal G_{-1})$, $\|S_n(f-m'')\|_2 \geq b_n\sqrt n$ infinitely many times.
\endproclaim

As shown by (e.g.) Wu and Zhao-Woodroofe, a fast rate of martingale approximation implies invariance priniciples and the law of
iterated logarithm (\cite{Wu07}, \cite{Z-Wo}). These limit theorems thus can be infered using the filtration $(\Cal F_j)$
but not (in the same way) using the natural filtration $(\Cal G_j)$.

\bigskip

\subheading{2. Proofs}

The proof of Proposition 1 leads to a problem of independent interest
which can be expressed as a property of analytic functions, or as a problem in spectral theory; these questions will be treated
in the last chapter of the article.

\demo{Proof of Proposition 1}
As shown in \cite{V93}, it is sufficient to prove that for every $m\in L^2(\Cal F_0) \ominus L^2 (\Cal F_{-1})$, 
$\limsup_{n\to\infty} \frac1{\sqrt n} \|\sum_{i=0}^{n-1}(X_i-m)\|_2 > 0$: if there is no $c>0$ in (1), there is a sequence
of $m_k\in L^2(\Cal F_0) \ominus L^2 (\Cal F_{-1})$ such that $\limsup_{n\to\infty} \frac1{\sqrt n} 
\|\sum_{i=0}^{n-1}(X_i-m_k)\|_2 \to 0$ as $k\to \infty$. Then, $m_k\to m\in L^2(\Cal F_0) \ominus L^2 (\Cal F_{-1})$
and $\lim_{n\to\infty} \frac1{\sqrt n} \|\sum_{i=0}^{n-1}(X_i-m_)\|_2 \to 0$.
\smallskip

Let $(e_i)_{i\in\Bbb Z}$ be a sequence of $\Cal N(0,1)$ distributed independent
random variables,
$a_i$, $i\in\Bbb Z$, be real numbers, $\sum_{i\in\Bbb Z} a_i^2 <\infty$, $X_n =
\sum_{i=0}^\infty a_ie_{n-i}$. $(X_n)$ is then a stationary linear process.
$T$ is a measure preserving and bimeasurable transformation such that $e_i
\circ T = e_{i+1}$.
$U$ is the unitary operator in $L^2(\Cal A)$ 
defined by $Uf= f \circ T$. We thus have $X_n = U^nX_0 = X_0\circ T^n$.
By $\Cal F_k$ we denote the filtration $\Cal F_k = \sigma\{e_i\,|\, i\leq k\}$.
Notice that $T^{-1}\Cal F_{k-1} = \Cal F_k$ and that the process $(X_k)$ is adapted to the filtration
$(\Cal F_k)$. We will find the process $(X_k)$ such that there is no martingale approximation w.r.t\.
$(\Cal F_k)$, and the random variables are mutually orthogonal. Because $(X_k)$ is a Gaussian process,
orthogonality implies independence.

\proclaim{Lemma 1 (\cite{V09})} If $m\in L^2(\Cal F_0)\ominus L^2(\Cal F_{-1})$ is such that
$$
  \frac1{\sqrt n} \|\sum_{i=0}^{n-1} (X_i-U^im)\|_2 \to 0
$$
then $m= ce_0$ for some $c\in\Bbb R$.
\endproclaim

\demo{Proof of Lemma 1} 
%The proof appears in \cite{V09}. 
We have $e_0\in L^2(\Cal F_0) \ominus L^2 (\Cal F_{-1})$; suppose that 
$E$ is the Hilbert space generated
by $e_0$ and $F$ its orthogonal complement in $L^2(\Cal F_0) \ominus L^2 (\Cal F_{-1})$, $m=m'+m''$ where
$m' = ce_0\in E$ for some real number $c$ and $m''\in F$. For $f=X_0$, the random variables $S_n(f-m')$ and $S_n(m'')$ are then 
mutually orthogonal, hence $\|S_n(m'')\|_2 \leq \|S_n(f-m)\|_2$, therefore $m''=0$.
\enddemo \qed

Without loss of generality we can suppose $\|X_0\|_2=1$. 

\proclaim{Lemma 2 (\cite{V09})} If $X_i$ are mutually orthogonal and
$$
  \limsup_{n\to\infty} \frac1{\sqrt n} \|\sum_{i=0}^{n-1}(X_i-ce_i)\|_2 =0,
$$
then $|c|=1$.
\endproclaim

\demo{Proof of Lemma 2} Recall that $\|X_0\|_2 =1 =\|e_0\|_2$. We have 
$\|\sum_{i=0}^{n-1}(X_i-ce_i)\|_2 \geq |\|\sum_{i=0}^{n-1} X_i\|_2 - \|\sum_{i=0}^{n-1} ce_i\|_2|
= |1-|c||\sqrt n$ by independence of $X_i$ and of $e_i$.
\enddemo \qed

\proclaim{Lemma 3}  Let $X_i$ be mutually orthogonal. Then 
$$
  \limsup_{n\to\infty} \frac1{\sqrt n} \|\sum_{i=0}^{n-1}(X_i-ce_i)\|_2 =0 \tag{1}
$$
if and only if 
$$
  \frac1{n} \sum_{k=0}^{n-1}\sum_{j=0}^{k} a_j \to c.  \tag{2}
$$
\endproclaim

\demo{Proof of Lemma 3} 
Recall that we assume $\|X_0\|_2 =1$. By Lemma 2, (1) can take place for $c=1$ or $c=-1$ only.
We have
$$
  \|\sum_{j=0}^{n-1} (X_j-ce_j)\|^2 = \|\sum_{j=0}^{n-1} X_j\|^2 + c^2 \|\sum_{j=0}^{n-1} e_j\|^2
  - 2c E[(\sum_{j=0}^{n-1} X_j)(\sum_{j=0}^{n-1} e_j)]; 
$$
because $(e_j)$ is an orthonormal sequence and $X_k = \sum_{j=0}^\infty a_je_{k-j}$, 
$$
  \limsup_{n\to\infty} \frac1{\sqrt n} \|\sum_{j=0}^{n-1}(X_j-ce_j)\|_2 = 0 
$$
if and only if 
$$
  \frac1{n} E[(\sum_{j=0}^{n-1} X_j)(\sum_{j=0}^{n-1} e_j)]] = 
  \frac1{n} \sum_{k=0}^{n-1}\sum_{j=0}^{k} a_j \to c.  \tag{2}
$$
\enddemo \qed

To finish the proof of Proposition 1 it is thus sufficient to find $(X_k)$ such that $X_k$ are 
mutually orthogonal, $\|X_0\|_2=1$, and $\frac1{n} \sum_{k=0}^{n-1}\sum_{j=0}^{k} a_j$ converge neither
to 1, nor to -1.
\medskip

Let $H$ be the Hilbert space generated by $e_i$, $i\in\Bbb Z$. Let
$K= \{z\in \Bbb C\,|\,|z|=1\}$ and $\lambda$ be the Lebesgue probability measure on $K$.
By $H'$ we denote the Hilbert space $L^2(\lambda)$; it is generated by the functions
$z^k$, $k\in\Bbb Z$. The multiplication $(Vf)(z) = zf(z)$ is a unitary operator in $H'$.
There exists a Hilbert space isomorphism $\phi$: $H\to H'$ such that $\phi(Uv)=V\phi(v)$ for
all $v\in H$ (\cite{A}). 
For any $f\in H'$ there thus exists $X=\phi^{-1}(f)\in H$; $X = \sum_{i\in\Bbb Z} a_{-i}e_i$
where $a_{-i} = \hat{f}(i)$ are the Fourier coefficients of $f$. $(X_k)=(U^kX)$ is thus a
gaussian process (\cite{A}). Because
$E(fU^kf) = \int_K f(z)\overline{z^kf(z)}\,d\lambda(z) = \int_K z^{-k}\,d\lambda(z)$,
the condition $|f|=1$ guarantees mutual orthogonality.

We thus seek a function $f(z) = \sum_{j=0}^\infty a_j z^{j}$ 
which is holomorphic on the open unit disc $D = \{z\in \Bbb C\,|\,|z|<1\}$, $a_j$ are all real, 
$\sum_{j=0}^\infty a_j^2 <\infty$, $|f|=1$ almost everywhere on $K$, 
and (2) holds neither for $c=1$, nor for $c=-1$. As shown in the next section, such a function exists.
\enddemo 
\hfill $\diamondsuit$
\medskip

\demo{Proof of Proposition 2}
Let $(e_i)_{i\in \Bbb Z}$ and $(f_i)_{i\in \Bbb Z}$ be two independent sequences of iid, each of them taking
values $\pm 1$ with probabilities $1/2$, $g = 1 + \sum_{i=1}^\infty e_{i}/3^{2i} + \sum_{i=1}^\infty e_{-i}/3^{2i+1}$,
$F = f_0g + 2e_0$.

We denote $\Cal E_k = \sigma\{e_i | i\leq k\}$, $\Cal F_k = \sigma\{f_i | i\leq k\}$, $\Cal C_k = \sigma\{e_i, f_i | i\leq k\}$,
and  $(\Cal G_i)$ is the natural filtration of the process $(F\circ T^i)$. Remark that the $\sigma$-algebra generated by $g$
equals $\Cal E_\infty$, hence $\Cal G_i = \Cal F_i \vee \Cal E_\infty$, the $\sigma$-algebra generated by $\Cal E_\infty$ and  
$\Cal F_i$. For any $i$, $g$ is $(\Cal G_i)$-measurable and we therefore have
$E(F|\Cal G_{-\infty}) = e_0$. Any $m\in L^2(\Cal G_0)\ominus L^2(\Cal G_{-1})$ is thus orthogonal w.r.t\. the sequence 
$(e_0\circ T^i)$ hence there does not exist any martingale approximation of the process $(F\circ T^i)$ w.r.t\. 
the natural filtration  $(\Cal G_i)$. On the other hand, there exists a nice approximation (a martingale-coboundary 
decomposition) for the filtration $(\Cal C_k)$. To see this, we notice that
\roster
\item"-" $E(f_0e_i|\Cal C_k) = f_ie_i$ if $k\geq 0$, $i\leq k$,
\item"-" $E(f_0e_i|\Cal C_k) = f_0E(e_i|\Cal C_k) =0$ if $k\geq 0$, $i> k$,
\item"-" $E(f_0e_i|\Cal C_k) = e_iE(f_0|\Cal C_k) =0$ if $k< 0$, $i\leq k$,
\item"-" $E(f_0e_i|\Cal C_k) = E(f_ie_i)=0$ if $k< 0$, $i> k$.
\endroster
We therefore have $E(F|\Cal C_{k+1}) - E(F|\Cal C_k) = f_0e_{k+1}/3^{2k+2}$ for all $k\geq 0$,
$E(F|\Cal C_0) - E(F|\Cal C_{-1}) = f_0e_0 + f_0 \sum_{i=1}^\infty e_{-i}/3^{2i+1} + 2e_0$,
$E(F|\Cal C_k) - E(F|\Cal C_{k-1}) = 0$ for all $k<0$. This implies that the Hannan's criterium 
$\sum_{k\in \Bbb Z} \|E(F|\Cal C_{k+1}) - E(F|\Cal C_k)\|_2 < \infty$ (cf\. \cite{Ha-He}; \cite{V93} for the non adapted version)
is satisfied and there is a martingale approximation.
\enddemo
\hfill $\diamondsuit$ 

\demo{Proof of Proposition 3} Recall that we denote $Uf = f\circ T$.
Let $e_k$, $1\leq k<\infty$, be random variables with $\|e_k\|_2 =1$, such that for each $k$, $(U^ie_k)_i$ is an iid sequence
and $(U^ie_k)_i$ are mutually independent processes. For each $k$, $e_k$ takes values $\pm q_k$ with probabilities $1/2q_k^2$
(the values of $q_k$ will be specified later).

Let $\rho_k$, $p_k$, $k=1,2,\dots$, be real numbers, $0<\rho_k \searrow 0$, $\rho_k<1/10$, $0<p_k$, $\sum_{k=1}^\infty p_k^2 < \infty$.
$0< \phi(k)\nearrow \infty$ are positive integers. We define
$$
  f = \sum_{k=1}^\infty p_k [\rho_ke_k - (1+\rho_k)U^{-\phi(k)}e_k].
$$
Notice that the series converges in $L^2$ (as we shall see later, the sequence of $q_n$ grows exponentially fast hence 
the series of $f$
converges almost surely). Suppose that the values of $q_k, p_k, \rho_k, \phi(k)$ are such that the 
random variables $e_k$ are 
$\sigma\{f\}$-measurable (this will be shown at the end of the proof). 
Therefore, $\Cal G_j = \sigma\{f\circ T^i : i\leq j\} =
\sigma\{e_k\circ T^i : i\leq j, k=1,2,\dots \}$; $(\Cal G_j)$ is the natural filtration of the process $(f\circ T^i)$.

Define $\Cal F_j = \sigma\{U^{i-\phi(k)}e_k : i\leq j, k=1,2,\dots \}$. $(\Cal F_j)$ is thus a filtration.

Define $m' = -\sum_{k=1}^\infty p_kU^{-\phi(k)}e_k \in L^2(\Cal F_0)\ominus L^2(\Cal F_{-1})$. By mutual orthogonality 
of the processes $(e_k\circ T^i)$ ($e_k\circ T^i$ are mutually independent and $E(e_k)=0$) and because $\|e_k\|_2=1$, we have
$$\gather
  \|S_n(f-m')\|_2^2 = \sum_{k=1}^\infty p_k^2 \rho_k^2 \|S_n(e_k -U^{-\phi(k)}e_k)\|_2^2 = \\
  = 2 \sum_{\phi(k)\leq n} p_k^2\rho_k^2\phi(k) + 2n \sum_{\phi(k)>n} p_k^2\rho_k^2.
  \endgather
$$

Define $m^* = -\sum_{k=1}^\infty p_ke_k \in L^2(\Cal G_0)\ominus L^2(\Cal G_{-1})$. Then
$$\gather
  \|S_n(f-m^*)\|_2^2 = \sum_{k=1}^\infty p_k^2(1+\rho_k)^2 \|S_n(e_k -U^{-\phi(k)}e_k)\|_2^2 = \\
  = 2 \sum_{\phi(k)\leq n} p_k^2(1+\rho_k)^2\phi(k) + 2n \sum_{\phi(k)>n} p_k^2(1+\rho_k)^2. 
  \endgather
$$

Because $\phi(k)\nearrow \infty$, $\sum_{k=1}^\infty p_k^2 < \infty$, and $\rho_k^2<1$,
$\|S_n(f-m')\|_2 = o(n)$ and $\|S_n(f-m^*)\|_2 = o(n)$.

Now, we will find the constants $p_k, \rho_k, \phi(k)$ so that $\|S_n(f-m'')\|_2 > b_n/\sqrt n$ infinitely many times for any 
$m''\in L^2(\Cal G_0)\ominus L^2(\Cal G_{-1})$, and $\|S_n(f-m')\|_2 \leq a_n$ for all $n$.

For the $p_k$ we can choose e.g\. $p_k = 1/k$. Then we define $\phi(0)=0$ and $\phi(k)\nearrow \infty$ increasing fast enough 
so that for every $j\geq 1$
$$
  2\sum_{k=j+1}^\infty p_k^2 > b_{\phi(j)}^2.
$$
This implies that $\|S_n(f-m^*)\|_2 \geq b_n\sqrt n$ for every $n=\phi(j)$. Notice that for any $m''\in 
L^2(\Cal G_0)\ominus L^2(\Cal G_{-1})$ we have $\|S_n(f-m'')\|_2 \geq \|S_n(m^*-m'')\|_2 - \|S_n(f-m^*)\|_2$ and 
$\|S_n(m^*-m'')\|_2 = \sqrt n \|m''-m^*\|_2$ hence, if $\|m''-m^*\|_2>0$, $\|S_n(f-m'')\|_2 > c\sqrt n$ infinitely many times 
for some $c>0$.

Next we find $\rho_k$ small enough so that $\|S_n(f-m')\|_2 \leq a_n$ for all $n\geq 1$. Without loss of generality we can 
suppose that $a_n=1$ for all $n$. 
We take $\rho_k^2 = 1/8\phi(k)$, for all $k$. Recall that $\sum_{k=1}^\infty 1/k^2 < 2$. Then 
$$
  2 \sum_{k=1}^\infty p_k^2\rho_k^2\phi(k) = \frac14 \sum_{k=1}^\infty  \frac1{k^2} < 1/2
$$
and
$$
  2n \sum_{\phi(k)>n} p_k^2\rho_k^2 < \frac14 \sum_{\phi(k)>n} \frac1{k^2} <1/2,\quad n=1,2,\dots ,
$$ 
hence $\|S_n(f-m')\|_2 \leq 1$. 

Eventually we choose the numbers $q_k$ so that the functions $e_k$ are $\sigma\{f\}$-measurable. 
Denote 
$$\gather
  f_{n} = \sum_{k=1}^n p_k [\rho_ke_k - (1+\rho_k)U^{-\phi(k)}e_k], \\
  r=r_n=3 \sum_{k=1}^n p_k q_k,\quad s=s_{n+1}=10r/\rho_{n+1} >100r.
  \endgather
$$ 
Notice that the previous calculation does not impose any condition on the values of $q_k$ and that
$|f_n| \leq r_n$.
We define the numbers $q_k$ so that $q_1 = 1$ and
$$
  10r_n = 30 \sum_{k=1}^n p_k q_k = \rho_{n+1}p_{n+1} q_{n+1},\quad n=1,2,\dots .
$$

We will show that for $x,y \in \{-q_{n+1}, 0, q_{n+1}\}$ the sets $\{e_{n+1}=x, U^{-\phi(n+1)}e_{n+1}=y\}$ are
$\sigma\{f_{n+1}\}$-measurable. By definition we get the inclusions

$\{e_{n+1}=-q_{n+1}, U^{-\phi(n+1)}e_{n+1}=-q_{n+1}\} \subset f_{n+1}^{-1}((s-r,s+r))$, 

$\{e_{n+1}=-q_{n+1}, U^{-\phi(n+1)}e_{n+1}=0\} \subset f_{n+1}^{-1}((-11r, -9r))$, 

$\{e_{n+1}=-q_{n+1}, U^{-\phi(n+1)}e_{n+1}=q_{n+1}\} \subset f_{n+1}^{-1}((-s-21r,-s-19r))$,

$\{e_{n+1}=0, U^{-\phi(n+1)}e_{n+1}=-q_{n+1}\} \subset f_{n+1}^{-1}((s+9r, s+11r))$,

$\{e_{n+1}=0, U^{-\phi(n+1)}e_{n+1}=0\} \subset f_{n+1}^{-1}((-r, r))$,

$\{e_{n+1}=0, U^{-\phi(n+1)}e_{n+1}=q_{n+1}\} \subset f_{n+1}^{-1}((-s-11r, -s-9r))$,

$\{e_{n+1}=q_{n+1}, U^{-\phi(n+1)}e_{n+1}=-q_{n+1}\} \subset f_{n+1}^{-1}((s+19r, s+21r))$,

$\{e_{n+1}=q_{n+1}, U^{-\phi(n+1)}e_{n+1}=0\} \subset f_{n+1}^{-1}((9r, 11r))$,

$\{e_{n+1}=q_{n+1}, U^{-\phi(n+1)}e_{n+1}=q_{n+1}\} \subset f_{n+1}^{-1}((-s-r, -s+r))$. \newline
Because the sets on the left give a partition of $\Omega$ and the sets on the right are mutually disjoint, the 
inclusions are equalities, hence the sets $\{e_{n+1}=x, U^{-\phi(n+1)}e_{n+1}=y\}$ are $\sigma\{f_{n+1}\}$-measurable. Therefore, 
the functions $e_{n+1}$ and $U^{-\phi(n+1)}e_{n+1}$ are $\sigma\{f_{n+1}\}$-measurable, $n\geq 1$. \newline
Recall $f_1 = p_1 [\rho_1e_1 - (1+\rho_1)U^{-\phi(1)}e_1]$ where $p_1 = 1 = q_1$, $\rho_1<1/10$, $e_1 =\pm 1$ with
probabilities $1/2$. We deduce that the functions $e_1$ and $U^{-\phi(1)}e_1$ are $\sigma\{f_1\}$-measurable.

Because the functions $e_{n+1}$ and $U^{-\phi(n+1)}e_{n+1}$ are $\sigma\{f_{n+1}\}$-measurable,
$f_n = f_{n+1} - p_{n+1} [\rho_{n+1}e_{n+1} - (1+\rho_{n+1})U^{-\phi({n+1})}e_{n+1}]$ is $\sigma\{f_{n+1}\}$-measurable.
%Because $e_n$ is $\sigma\{f_{n}\}$-measurable 
%(for $n\geq 2$ this follows from the preceding calculation and for $n=1$ we deduce this 
%from the definition of $f_n$, $q_1 = 1$ and $\rho_1<1/10$), 
Using $\sigma\{f_{n}\}$-measurability of $e_n$ and $U^{-\phi(n)}e_{n}$ we deduce that
$e_n$ and $U^{-\phi(n)}e_{n}$ are $\sigma\{f_{n+1}\}$-measurable. 
By induction we deduce that all $e_k$ (and $U^{-\phi(k)}e_{k}$ too) with $1\leq k\leq n+1$ are $\sigma\{f_{n+1}\}$-measurable.

The $q_n$ are growing exponentially fast hence the measures of the supports of $e_n$ are decreasing exponentially fast. For
$F_n = \cup_{k=n}^\infty \{e_k \neq 0\}$ we thus have $\mu(F_n)\searrow 0$. For $E_n = \Omega\setminus F_n$
we have $E_n\subset E_{n+1}$ for all $n$ and $\mu(E_n)\nearrow 1$. Because $f_n=f$ on $E_n$, $f_n^{-1}(B) \cap E_n =
f^{-1}(B) \cap E_n$ for every Borel set $B\subset \Bbb R$. Therefore, $\sigma\{e_k\}\cap E_n \subset 
\sigma\{f\}\cap E_n$, $1\leq k\leq n$, $n=1,2,\dots$ . We deduce that all the functions $e_k$, $k\geq 1$, are 
$\sigma\{f\}$-measurable.

\enddemo
\hfill $\diamondsuit$

\bigskip
\subheading{3. A bad approximation property}

In the preceding section, the problem of finding a gaussian, adaptated process  $(X_k)$, which is in some sense (described before) 
badly approximable, was reduced to the following, function-theoretical, question, in which one sets 

$$f(z)=\sum_{n=0}^\infty a_{n}z^n,\ A_n=\sum_{j=0}^{n}a_j, \ M_n={A_0+\ldots+A_{n-1}\over n}, f^{*}(e^{it})=
\lim_{r\buildrel <\over{\to 1}}f(re^{it}),$$

the function $f$  being a bounded, analytic function in the open unit disk $D$ of the complex plane, with {\it real} coefficients 
$a_n$, and $f^*(e^{it})$ denoting its radial limits, which exists almost everywhere with respect to the Lebesgue measure, according 
to a well-known result (Fatou's theorem: see \cite{DU} p.6 or \cite{RU} p.340). 
If this radial limit has modulus one (a.e.), the function $f$ is said to be {\it inner} (see \cite{DU} p.24 or \cite{RU} p.342), 
which is a severe restriction on its behavior. But this is this "innerness" property which has guaranteed us  the orthogonality, 
therefore the mutual  independence, of the $X_k's$.  And the problem now is :

 {\bf Problem} Find an inner function $f$ with real coefficients  such that one of the following cases occurs : 
\roster
  \item $M_n$ has no limit,
  \item $M_n$ has a limit $c$, but $c\neq\pm 1$.
\endroster
  This turns out to be  possible in many ways. We just  indicate two examples below, of a different nature.
 \bigskip

\bigskip

 {\bf Example 1.} Let $a>0$ and $f(z)=e^{-a({1+z\over 1-z})}$. This is clearly a zero-free inner function (called singular) with real coefficients $a_n$, and with 
$f^{*}(e^{it}) =e^{-ia {\ cotg}{t\over 2} }$. Now, it was proved in \cite{NESH} that we have 
$$a_n=\pi^{-{1\over 2}}(2a)^{{1\over 4}}n^{-{3\over 4}}{\ cos}\Big( 2\sqrt{2an}+{\pi\over 4}\Big)+O(n^{-{5\over 4}}),$$ from which it 
easily follows, using summation by parts, that $A_n$ converges, necessarily to $0$ by Abel's theorem, since the radial limit of $f$ at 
$1$ is clearly $0$. (Another approach to that example is given in \cite{BAKOQU} where it is observed that $A_n=e^{-a}L_{n}(2a)$, where 
$L_n$ is the $n$-th Laguerre polynomial, so that $A_n\to 0$).
Therefore, $M_n\to 0$, i.e. $M_n$ has a limit $c$, but this limit does not have the correct value $c=\pm 1$, so that the corresponding 
gaussian process is badly approximable.
 \bigskip
 
 {\bf Example 2.} Let $(z_n)_{n\geq 1}$ be a sequence of real numbers between $0$ and $1$, such that $\sum_{n=1}^\infty (1-z_n)<\infty$ 
(a Blaschke sequence) and let  $B$ be the corresponding Blaschke product, namely 
 $$B(z)=\prod_{n=1}^\infty {z_n -z\over 1-z_{n}z}.\eqno (1)$$
 It is well-known(\cite{RU} p.312) that this is an inner function. Assume that $M_n$ has a limit $c$. Then, by a well-known extension 
(due to Frobenius) of Abel's theorem (\cite{KOR} p.4), we have $\lim_{r\buildrel <\over{\to}1} B(r)=c$. This last fact may quite well 
happen, taking for example  $z_n=1-n^{-\alpha}, \alpha>1$ (\cite{RU} p.317). But then, necessarily, $c=0$ since $B(z_n)=0$ and $z_n\to 1$. Observe that a finite Blaschke product, even with complex zeros symmetric with respect to the real axis, would never do the job, since then $B$ extends analytically across the closed unit disk, and $A_n\to B(1)=1$, so that $M_n\to 1$. (This might be the only case of an inner function with real coefficients for which 
 $M_n\to 1$). Finally, if $M_n$ has no limit, we automatically  have an example of a badly approximable gaussian process. We show that this can happen on a simple specific example:
 \medskip

\proclaim{Proposition 6} There exists an infinite  Blaschke product $B(z)=\sum_{n=0}^\infty a_n z^n$ with real zeros and coefficients,  
verifying the two following conditions : 
\roster

\item $ { \lim\inf}_{r\buildrel<\over{\to 1}}\vert B(r)\vert =0\,\, \text{and}\ { \lim\sup}_{r\buildrel<\over{\to 1}}\vert B(r)\vert >0$;
\item $A_n  =a_0+\ldots +a_n$ is not Ces\`aro-summable, i.e. $M_n$ has no limit.
\endroster
\endproclaim

\demo{Proof} Observe that 2. is an automatic consequence of 1., through the Frobenius theorem already mentioned.  We now take $$z_n=1-2^{-n}.$$  Let $r={z_n +z_{n+1}\over 2}$, and proceed to minorize $B(r)$. We have 
$\vert B(r)\vert =P_1 P_2P_3P_4$, where 
$$P_1=\prod_{j<n} {r-z_j\over 1-z_{j}r};\ P_2= {r-z_n\over 1-z_{n}r};P_3= {-r+z_{n+1}\over 1-z_{n+1}r};
P_4=\prod_{j>n+1} {-r+z_j\over 1-z_{j}r}.$$
Now, we see that 
\roster
\item   $$P_1\geq \prod_{j<n} {z_n-z_j\over 1-z_{j}z_n}\geq \prod_{j<n}{2^{-j}-2^{-n}\over 2^{-j}+2^{-n}}= 
\prod_{k<n}{1-2^{-k}\over 1+2^{-k}}$$
$$\geq \prod_{k=1}^\infty{1-2^{-k}\over 1+2^{-k}}=c>0.$$
\item One has similarly $P_4\geq c>0$.
\item $$P_2= {z_{n+1}-z_n\over 2-z_n(z_n+z_{n+1})}\geq {z_{n+1}-z_n\over 2(1-z_{n}^2)}\geq {1\over 4}{z_{n+1}-z_n\over (1-z_{n})}={1\over 8}.$$
\item One has similarly $P_3\geq {1\over 8}.$
\endroster
This ends the proof of Proposition 6.
\enddemo

\underbar{Remarks :} 
\roster
\item  The specific sequence $z_n=1-2^{-n}$  in the proof of the proposition verifies  the Newman lacunarity condition: 
$${1-\vert z_{n+1}\vert \over 1-\vert z_n\vert}\leq \rho<1.$$ Such sequences will in particular be interpolation sequences in 
the sense of 
Carleson (\cite{DU} p.149 or \cite{HOF} p.203), and the Taylor coefficients of  $a_n$ of the corresponding Newman-lacunary 
Blaschke product will verify (\cite{NESH}):

% \begin{equation}\label{Small} a_n=O({1\over n}).\end{equation}

$$
  a_n=O({1\over n}).
$$

\item It can be proved that, for any infinite Blaschke product, one has $$\overline{\lim_{n\to\infty}}n\vert a_n\vert>0.$$ 
(see \cite{NESH}). Therefore, the estimate $a_n=O({1\over n})$ is the best possible.

\item The estimate on the Taylor coefficients of a  Blaschke product with Newman-lacunary zeros was recently used by Bourgain 
and Kahane (\cite{BOKA}). 

\endroster

\enddocument                                                                    
\end